\documentclass{amsart}

\usepackage{amsfonts}
\usepackage{amsmath}
\usepackage{amsthm}
\usepackage{amssymb}
\usepackage[english]{babel}
\usepackage{graphics}
\usepackage{latexsym}
\usepackage{longtable} 
\usepackage{mathrsfs}
\usepackage{graphicx}
\usepackage{wrapfig}
\usepackage[pdftex,colorlinks=true,linkcolor=blue,citecolor=cyan]{hyperref}
\usepackage{ upgreek }

\usepackage[utf8x]{inputenc}
\usepackage{amsmath}
\usepackage{tikz}

\usepackage{enumerate}

\newtheorem{theorem}{Theorem}[section]
\newtheorem{corollary}[theorem]{Corollary}
\newtheorem{lemma}[theorem]{Lemma}
\newtheorem{proposition}[theorem]{Proposition}

\newtheorem{claim}[theorem]{Claim}
\newtheorem{example}[theorem]{Example}

\theoremstyle{definition}
\newtheorem{definition}[theorem]{Definition}

\newcommand{\mD}{\mathcal{D}}

\newcommand{\R}{\mathbb{R}}
\newcommand{\N}{\mathbb{N}}
\newcommand{\mS}{\mathbb{S}}
\newcommand{\mB}{\mathbb{B}}

\newcommand{\D}{\mathrm{D}}

\newcommand{\ms}{\medskip}

\newcommand{\al}{\alpha}
\newcommand{\be}{\beta}
\newcommand{\ga}{\gamma}
\newcommand{\de}{\delta}

\newcommand{\Om}{\Omega}

\newcommand{\larrow}{\longrightarrow}
\newcommand{\ot}{\otimes}

\newcommand{\p}{\partial}
\newcommand{\sub}{\subseteq}
\newcommand{\set}{\setminus}
\newcommand{\by}{\times}
\newcommand{\rk}{\mathrm{rk}}

\newcommand{\Div}{\mathrm{Div}}

\newcommand{\bt}{\begin{theorem}}\newcommand{\et}{\end{theorem}}
\newcommand{\bd}{\begin{definition}}\newcommand{\ed}{\end{definition}}
\newcommand{\bl}{\begin{lemma}}\newcommand{\el}{\end{lemma}}
\newcommand{\beq}{\begin{equation}}\newcommand{\eeq}{\end{equation}}
\newcommand{\bc}{\begin{claim}}\newcommand{\ec}{\end{claim}}
\newcommand{\bex}{\begin{example}}\newcommand{\eex}{\end{example}}
\newcommand{\bcor}{\begin{corollary}}\newcommand{\ecor}{\end{corollary}}
\newcommand{\bp}{\begin{proof}}\newcommand{\ep}{\end{proof}}

\newcommand{\BPT}{\medskip \noindent \textbf{Proof of Theorem} }

\numberwithin{equation}{section}

\begin{document}

\title[Phase separation of $ \textbf{n-}$dimensional  $\infty$ \!\!-\!\! Harmonic  mappings]{Phase separation of $ \textbf{n-}$dimensional  $\infty$ \!\!-\!\! Harmonic  mappings}

\author{Hussien Abugirda}

\address{Department of Mathematics , College of Science, University of Basra, Basra, Iraq AND Department of Mathematics and Statistics, University of Reading, Whiteknights, \\ PO Box 220, Reading RG6 6AX, UK}

\email{h.a.h.abugirda@pgr.reading.ac.uk}


\begin{abstract} Among other interesting results, in a  recent paper \cite{15}, Katzourakis analysed the phenomenon of separation of the solutions $ u\!:\! \R^2 \supseteq \Om \longrightarrow \R^{N}$, to the $\infty$-\!Laplace system
\beq 
\Delta_\infty u \, :=\, \Big(\D u \ot \D u+|\D u|^{2} [\![\D u]\!]^\bot \ot I\Big):\D^{2}u\,=\,0 ,  \nonumber
\eeq
to phases with qualitatively different behavior in the case of $n=2 \leq N$. The solutions of the $\infty$-\!Laplace system are called the  $\infty$ \!\!-\!\! Harmonic  mappings.
In this paper we discuss an extension of Katzourakis’ result  mentioned above to higher dimensions by studying the phase separation of $n$-dimensional $\infty$ \!\!-\!\! Harmonic  mappings in the case $N \geq n \geq 2$.
\end{abstract}

\maketitle

\section{Introduction} \label{section1}

In this paper we study the phase separation of $n$-dimensional  $\infty$ \!\!-\!\! Harmonic  mappings  $ u\!:\! \R^n \supseteq \Om \longrightarrow \R^{N}$ by which we mean the classical solutions to the $\infty$-\!Laplace system
\beq \label{1.1}
\Delta_\infty u \, :=\, \Big(\D u \ot \D u+|\D u|^{2} [\![\D u]\!]^\bot \ot I\Big):\D^{2}u\,=\,0 , \,\,\,\,\,\, \text{on} \ \Om, 
\eeq
where  $n,N$ are integers such that $N \geq n \geq 2$ and $\Om$ an open subset of $ \R^n$.
Here, for the map $u$ with components $(u_1,...,u_N)^\top$ the notation $\D u$ symbolises the gradient matrix
\beq \label{1.2}
\ \ \ \ \D u(x) = \big(\D_i u_\al(x)\big)_{i=1...n}^{\al=1...N} \, \in\, \R^{N\by n}\ , \ \ \D_i \equiv \p /\p x_i,
\eeq
and for any $X\in \R^{N\by n}$, $ [\![X]\!]^\bot$ denotes the orthogonal projection on the orthogonal complement of the range of linear map $X :\R^n \larrow \R^N$:
\beq \label{1.3}
[\![X]\!]^\bot := \textrm{Proj}_{\mathrm{R}(X)^\bot}.
\eeq
In index form, the system \eqref{1.1} reads
\[
\ \ \ \sum_{\be=1}^N\sum_{i,j=1}^n \Big(\D_i u_\al \, \D_j u_\be + |\D u|^2[\![\D u]\!]_{\al\be}^\bot\, \de_{ij}\Big)\D^2_{ij}u_\be\,=\,0, \  \ \ \al=1,...,N.
\]

Our general notation will be either self-explanatory, or otherwise standard as e.g.\ in \cite{9,10,28}. Throughout this paper we reserve $n, N \in \N$ for the dimensions of Euclidean spaces and $\mathbb{S} ^{N-1}$ denotes the unit sphere of $\R ^N$.

Speaking about the system \eqref{1.1}, we would like to mention that the system  \eqref{1.1} is called the``$\infty$-\!Laplacian" and it arises as a sort of Euler-Lagrange PDE of vectorial variational problems in $L^\infty$ for the supremal functional
\beq \label{1.4}
\ \ \ E_\infty(u,\mathcal{O}) :=  \|H(\D u)\|_{L^\infty(\mathcal{O})}, \ \ \ u \in W_{loc}^{1,\infty}(\Om,\mathbb{R}^N), \ \mathcal{O} \Subset \Omega,
\eeq
when the Hamiltonian (the non-negative function $H \in  C^{2}\!\left( \R^{N \by n}\right) $) is chosen to be $H(\D u) \,=\,\frac{1}{2} \big|\D u \big|^2 $, with $\big| . \big|$ is the Euclidean norm on the space $ \R^{N \by n}$.  the $\infty$-\!Laplacian is a special case of the system 
\beq \label{1.5}
\Delta_\infty u \, :=\, \Big(H_P \ot H_P+ H [\![H_P]\!]^\bot H_{PP}\Big)(\D u):\D^{2}u\,=\,0 , 
\eeq
which was first formally derived by Katzourakis  \cite{11} as the limit of the Euler-Lagrange equations of the integral functionals $E_m (u,\Om)\, :=\, \int_\Om(H(\D u))^p$ as $p \larrow \infty$.

Eventhough the theory of weak solutions has witnessed a significant development so far, particularly the new theory of ``$\mathcal{D}$-solutions"  introduced by Katzourakis  \cite{20}, which applies to nonlinear PDE systems of any order and allows for merely measurable maps to be rigorously interpreted and studied as solutions of PDE systems fully nonlinear and with discontinuous coefficients,  yet the structure of weak solutions are complicated to understand. In this paper, we restrict our attention to classical solutions which might be helpful to imagine and understand the behavior and the structure of the weak solutions.

For the $\infty$-\!Laplace system \eqref{1.1}  the orthogonal projection on the orthogonal complement of the range,  $ [\![\D u]\!]^\bot$, coincides with the projection on the geometric normal space of the image of the solution.

It is worth noting that $\infty$ \!\!-\!\! Harmonic  maps are affine when $n=1$ since in this case the system \eqref{1.1} simplifies to 
\beq \label{1.6}
\Delta_\infty u \, =\, \Big( u' \ot u' \Big) u''+|u'|^{2}  \Big (I - \frac{ u' \ot u'}{|u'|^{2}}\Big) u'' \,=\, |u'|^{2}u'' , 
\eeq
and hence no interesting phenomena arise when $n=1$.

For the case $N=1$, the system \eqref{1.1} reduced to the  single $\infty$-\!Laplacian PDE 
\beq \label{1.7}
\Delta_\infty u \, :=\, \Big(\D u \ot \D u\Big):\D^{2}u\,=\,0,
\eeq
 since the normal coefficient $|\D u|^{2} [\![\D u]\!]^\bot $ vanishes identically. This also happen when $u$ is submersion. 
The single $\infty$-\!Laplacian PDE \eqref{1.7}, and the related scalar  $L^\infty$-\!variational problems, started being studied in the '60s  by Aronsson in  \cite{4,5}. Today it is being studied in the context of Viscosity Solutions  (see for example Crandall  \cite{3}, Barron-Evans-Jensen  \cite{7} and Katzourakis \cite{16}).

 The vectorial case $N  \geq 2$ first arose in the early 2010s  in the work of Katzourakis \cite{11}. Due to both the mathematical significance as well as the importance for applications particularly in Data Assimilation, the area is developing very rapidly (see \cite{2}, \cite{6},\cite{8}, \cite{12}-\cite{15}, and also \cite{17}-\cite{27}).

In a joint work with Katzourakis and Ayanbayev \cite{1}, among other results, we have proved that the image $u(\Om)$ of a solution $u \in C^2(\Om,\R^N)$  to the nonlinear system \eqref{1.1}  satisfying that the rank of its gradient matrix is at most one, $\rk(\D u) \leq 1 \, \text{ in }\Om$, is contained in a polygonal line in $\R^N$, consisting of an at most countable union of affine straight line segments (possibly with self-intersections). Hence the component  $ [\![\D u]\!]^\bot \Delta u$ of $\Delta_\infty$ forces flatness of the image of solutions.

 Interestingly, even when the operator $\Delta_\infty$ is applied to $C^{\infty}$ maps, which may even be solutions, \eqref{1.1} may have discontinuous coefficients. This further difficulty of the vectorial case is not present in the scalar case.
 As an example consider 
\beq \label{1.8}
u(x,y):=e^{ix}-e^{iy} ,   \ \ \ \ u: \R^2 \larrow \R^2.
\eeq

 Katzourakis has showed in  \cite{11} that even though  \eqref{1.8} is a smooth solution of the  $\infty$-\!Laplacian near the origin, we still have the coefficient $|\D u|^{2} [\![\D u]\!]^\bot $ of \eqref{1.1} is discontinuous. This is because  when the dimension of the image changes, the projection $[\![\D u]\!]^\bot$ ``jumps''. More precisely, for  \eqref{1.8}  the domain  splits to three components according to the $\rk(\D u)$, the ``2D phase $\Om_2$'', whereon $u$ is essentially 2D, the``1D phase $\Om_1$'', whereon $u$ is essentially $1D$ (which is empty for \eqref{1.8})and the ``interface $S$'' where the coefficients of $\Delta_\infty$ become discontinuous. 

 In \cite{12} Katzourakis constructed  additional examples, which are  more intricate than \eqref{1.8}, namely smooth 2D  $\infty$ \!\!-\!\! Harmonic  maps whose interfaces have triple junctions and corners and are given by the explicit formula 
\beq \label{1.9}
u(x,y):= \int_{y}^{x} e^{iK(t)}dt.  
\eeq
 Indeed, for $K \in  C^{1}\!\left( \R, \R\right)$ with $\| K\|_{L^\infty(\R)} < \frac{\pi}{2}$, \eqref{1.9} defines $ C^{2}$  $\infty$ \!\!-\!\! Harmonic  map whose phases are as shown in Figures 1(a), 1(b) below, when $K$ qualitatively behaves as shown in the Figures 2(a), 2(b) respectively \footnote{The figures 1(a), 1(b),  2(a) and 2(b) are courtesy of N. Katzourakis.}. Also, on the phase $\Om_1$ the $\infty$ \!\!-\!\! Harmonic  map  \eqref{1.9} is given by a scalar $\infty$ \!\!-\!\! Harmonic function times a constant vector, and on the phase $\Om_2$ it is a solution of the vectorial Eikonal equation. 
\begin{align}
& \underset{\scriptstyle{\text{Figure 1(a). \hspace{100pt} Figure 1(b).\ \ \  }}}{ \includegraphics[scale=.18]{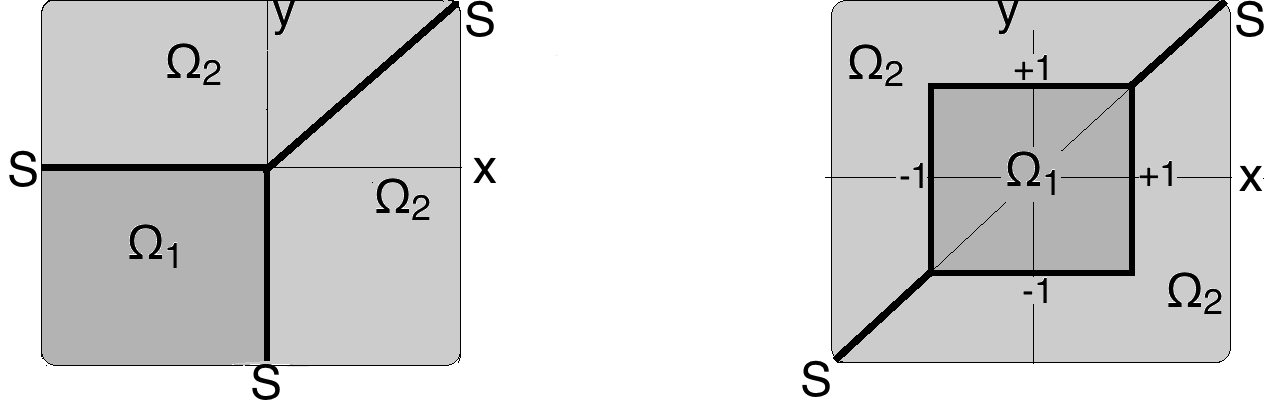}} \nonumber\\
& \underset{\scriptstyle{\text{Figure 2(a). \hspace{100pt} Figure 2(b).}}}{  \includegraphics[scale=.16]{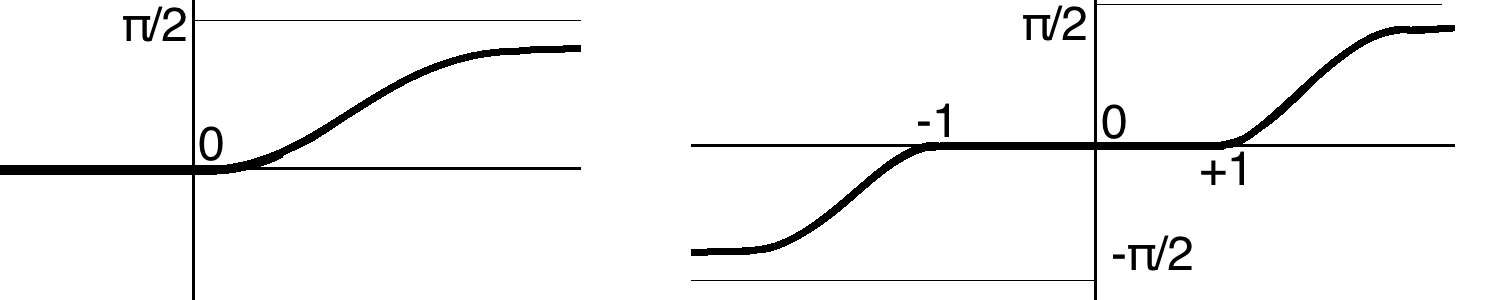} } \nonumber
\end{align}

One of the interesting results in  \cite{15} was that this phase separation is a general phenomena for smooth 2D  $\infty$ \!\!-\!\! Harmonic  maps. Therein the author proves that on each phase the dimension of the tangent space is constant and these phases are separated by interfaces whereon $[\![\D u]\!]^\bot$ becomes discontinuous. Accordingly the author established the next result:

\begin{theorem} [Structure of 2D $\infty$ \!\!-\!\! Harmonic  maps, cf.\ \cite{15}] \label{theorem1.1}\

Let $ u\!:\!  \R^2 \supseteq  \Om \longrightarrow \R^{N}$  be an $\infty$ \!\!-\!\! Harmonic  map in  $ C^{2}\!\left( \Om,\R^N\right)$, that is a solution to \eqref{1.1}. Let also $N \geq 2$. Then, there exist disjoint open sets $\Om_{1}$, $\Om_{2} \sub \Om$, and a closed nowhere dense set $S$ such that $\Om\,=\, \Om_{1} \bigcup S \bigcup \Om_{2}$ and:
\begin{enumerate} [(i)] 
\item On $\Om_{2}$ we have  $\rk(\D u)\,=\,2\ $, and the map $ u\!:\! \Om_{2} \longrightarrow \R^{N}$ is an immersion and solution of the Eikonal equation:
\beq \label{1.10}
|\D u|^{2}= C^{2} \textgreater 0.
\eeq
The constant $C$ may vary on different connected components of $\Om_{2}$.
\item  On $\Om_{1}$ we have $\rk(\D u)\,=\,1$ and the map $ u\!:\! \Om_{1} \longrightarrow \R^{N}$ is given by an essentially scalar $\infty$ \!\!-\!Harmonic function $f\!:\! \Om_{1} \longrightarrow \R$:
\beq \label{1.11}
u\,=\,a+\xi f, \,\,\, \Delta_\infty f \,=\,0, \,\,\,      a \in \R^{N},\,\,\, \xi \in  \mathbb{S} ^{N-1}.
\eeq
The vectors $a,\xi$ may vary on different connected components of $\Om_{1}$.
\item  On S, $|\D u|^{2}$ is constant and also $\rk(\D u)=1$. Moreover if $S= \partial \Om_1 \cap \partial \Om_2$ (that is if both the $1\D$ and $2\D$ phases coexist) then  $ u\!:\! S \longrightarrow \R^{N}$ is given by an essentially scalar solution of the Eikonal equation:
\beq \label{1.12}
u\,=\,a+\xi f, \,\,\, |\D f|^{2}= C^{2} \textgreater 0, \,\,\,      a \in \R^{N},\,\,\, \xi \in  \mathbb{S} ^{N-1}.
\eeq
\end{enumerate}

\end{theorem}

\ms

\ms

 The main result of this paper is to generalise these results to higher dimension $N \geq n \geq 2$. The principle result in this paper in the following extension of theorem \ref{theorem1.1}:

\begin{theorem}[ Phase separation of $n$-dimensional  $\infty$ \!\!-\!\! Harmonic  mappings] \label{theorem1.2}\

Let $\Om\sub \R^n$ be a bounded open set, and let $ u\!:\! \Om \longrightarrow \R^{N}$, $N \geq n \geq 2$, be an $\infty$ \!\!-\!\! Harmonic  map in $ C^{2}\!\left( \Om,\R^N\right)$, that is a solution to the $\infty$-\!Laplace system  \eqref{1.1}. Then, there exist disjoint open sets $ \big ( \Om_{r} \big ) _{r=1}^{n} \sub \Om$, and a closed nowhere dense set $S$ such that $\Om\,=\, S \bigcup \Big( \bigcup \limits_{i=1}^{n} \Om_{i} \Big)$ such that:
\begin{enumerate}[(i)] 
\item On $\Om_{n}$ we have $\rk(\D u)\,\equiv \,n$ and the map $ u\!:\! \Om_{n} \longrightarrow \R^{N}$ is an immersion and solution of the Eikonal equation:
\beq \label{1.13}
|\D u|^{2}= C^{2} \textgreater 0.
\eeq
The constant $C$ may vary on different connected components of $\Om_{n}$.
\item  On $\Om_{r}$ we have $\rk(\D u)\,\equiv\,r$,  where $r$ is integer in $\{2,3,4,...,(n-1) \}$, and $|\D u(\gamma(t))|$ is constant along trajectories of the parametric gradient flow of $u( \ga(t,x,\xi))$ 
\beq \label{1.14} 
\left\{
\begin{array}{ll}
\dot \ga(t,x,\xi)= \xi ^ {\top}  \D u \big(\ga(t,x,\xi) \big), &   t\in (-\upvarepsilon,0) \bigcup (0,\upvarepsilon), \ms\\
\ga (0,x,\xi) = x,
\end{array}
\right.
\eeq
where $\xi \in  \mathbb{S} ^{N-1}$, and  $\xi \notin N\big(\D u \big(\ga(t,x,\xi) \big)^\top \big)$.
\item  On $\Om_{1}$ we have $\rk(\D u)\,\leq\,1$ and the map $ u\!:\! \Om_{1} \longrightarrow \R^{N}$ is given by an essentially scalar $\infty$ \!\!-\!Harmonic function $f\!:\! \Om_{1} \longrightarrow \R$:
\beq \label{1.15}
u\,=\,a+\xi f, \,\,\, \Delta_\infty f \,=\,0, \,\,\,      a \in \R^{N},\,\,\, \xi \in  \mathbb{S} ^{N-1}.
\eeq
The vectors $a,\xi$ may vary on different connected components of $\Om_{1}$.
\item  On S, when $S \supseteq \partial \Om_p \cap \partial \Om_q \,=\, \phi$ for all $p$ and $q$ such that $2 \leq p < q \leq n-1$, then we have that $|\D u|^{2}$ is constant and also $\rk(\D u) \,\equiv\, 1$. Moreover on 
\[
 \partial \Om_1 \cap \partial \Om_n \sub S,
\]
(when both $1\D$ and $n\D$ phases coexist), we have that  $ u\!:\! S \longrightarrow \R^{N}$ is given by an essentially scalar solution of the Eikonal equation:
\beq \label{1.16}
u\,=\,a+\xi f, \,\,\, |\D f|^{2}= C^{2} \textgreater 0, \,\,\,      a \in \R^{N},\,\,\, \xi \in  \mathbb{S} ^{N-1}.
\eeq
On the other hand, if there exist some $r$ and $q$ such that $2 \leq r < q \leq n-1$, then on $S \supseteq \partial \Om_r \cap \partial \Om_q \, \neq\, \phi$ (when  both $r\D$ and $q\D$ phases coexist), we have that $\rk(\D u)\,\equiv\,r$ and we have same result as in (ii) above.
\end{enumerate}

\end{theorem}

\ms
\ms

\section{Preliminaries} \label{section2}
For the  convenience of the reader, in this section we recall without proof a theorem of rigidity of rank-one maps, proved in \cite{15}, which will be used in the proof of the main result of this paper in section \ref{section3}. We also recall the proposition of Gradient flows for tangentially $\infty$ \!\!-\!\! Harmonic maps which introduced in \cite{11} and its  improved modification lemma in \cite{14}.
\begin{theorem} [Rigidity of Rank-One maps, cf.\ \cite{15}] \label{theorem2.1}\

Suppose $\Om\sub \R^n$ is open and contractible and  $ u\!:\! \Om \longrightarrow \R^{N}$ is in   $ C^{2}\!\left( \Om,\R^N\right)$. Then the following are equivalent:
\begin{enumerate}[(i)] 
\item $u$ is a Rank-One map, that is $rk(\D u)\,\leq \,1$ on $\Om$ or equivalently there exist maps  $ \xi\!:\! \Om \longrightarrow \R^{N}$ and  $ w\!:\! \Om \longrightarrow \R^{n}$ with $w \in  C^{1}\!\left( \Om,\R^n\right)$ and $\xi \in  C^{1}\!\left( \Om \set \{w=0\},\R^N\right)$ such that $\D u \,=\, \xi \ot w$.
\item There exist $f \in  C^{2}\!\left( \Om,\R\right)$, a partition ${\big\{ B_i \big\} }_{i \in \N}$ og $\Om$ to Borel sets where each $ B_i$ equals a connected open set with a boundary portion and Lipschitz curves  ${\big\{\mathcal{V}^i \big\} }_{i \in \N} \subseteq W^{1,\infty}_{\text{loc}}(\Om)^N$ such that on each $ B_i $ u equals the composition of $\mathcal{V}^i$ with $f$:
\beq \label{2.1}
u\,=\, \mathcal{V}^i  \circ f  \ \ \ , \ \ \ \ \  \text{on} \,  B_i  \subseteq \Om.
\eeq
Moreover, $|\dot {\mathcal{V}}^i| \, \equiv \, 1$ on $f(B_i)$, $\dot {\mathcal{V}}^i \, \equiv \, 0$ on $\R \set f(B_i)$ and there exist $\ddot  {\mathcal{V}}^i$ on $ f(B_i)$, interpreted as 1-sided on $\partial f(B_i)$, if any. Also,
\beq \label{2.2}
\D u \, = \, ( \mathcal{V}^i  \circ f )   \ot \D f  \ \ \ , \ \ \ \ \  \text{on} \,  B_i  \subseteq \Om,
\eeq
and the image $u(\Om)$ is an 1-rectifiable subset of $\R^N$:
\beq \label{2.3}
u(\Om) \, = \,  \bigcup \limits_{i=1}^{\infty}  \mathcal{V}^i (f(B_i)) \subseteq \R^N.
\eeq
\end{enumerate}

\end{theorem}

\ms

\ms

\begin{proposition} [Gradient flows for tangentially $\infty$ \!\!-\!\! Harmonic maps, cf.\ \cite{11}] \label{theorem2.2}\

Let  $ u \in C^{2}\!\left( \R^n,\R^N\right)$. Then, $\D u \, \D \Big( \frac{1}{2}|\D u|^{2} \Big)\,=\,0$ on $\Om \Subset \R^n $ if and only if the flow map $ \ga\!:\! \R \times \Xi \longrightarrow \Om$ with $\Xi := \{ (x,\xi) \ | \ \xi^{\top} \D u(x) \neq 0 \} \sub \Om \times \mathbb{S} ^{N-1}$ of 
\beq \label{2.4} 
\left\{
\begin{array}{ll}
\dot \ga(t,x,\xi)= \xi ^ {\top}  \D u \big(\ga(t,x,\xi) \big), \ms\\
\ga (0,x,\xi) = x,
\end{array}
\right.
\eeq
satisfies along trajectories
\beq \label{2.5} 
\left\{
\begin{array}{ll}
\big|\D u \big( \ga(t,x,\xi) \big)\big|= \big|\D u \big( \ga(x) \big)\big|,\ \ \ t \in \R \ms\\
t \longmapsto \xi ^ {\top}  u \big(\ga(t,x,\xi) \big) \ \text{ is increasing}.
\end{array}
\right.
\eeq
\end{proposition}

\ms

\ms

The following lemma is improved modification of proposition  \ref{theorem2.2}  

\bl [cf.\ \cite{14}] \label{lemma2.3}\

Let $ u\!:\!  \R^n \supseteq  \Om \longrightarrow \R^{N}$ be in  $ u \in C^{2}\!\left( \Om,\R^N\right)$. Consider the gradient flow
\beq \label{2.6} 
\left\{
\begin{array}{ll}
\dot \ga(t,x,\xi)=\Big( \frac{|\D u|^{2}}{| \xi ^ {\top}\D u|^{2}} \xi ^ {\top}  \D u \Big) \big(\ga(t,x,\xi) \big), \ \ \ t \neq 0 \ms\\
\ga (0,x,\xi) = x,
\end{array}
\right.
\eeq
for $x \in \Om$, $\xi \in \mathbb{S} ^{N-1} \set [\![\D u]\!]^\bot$. Then, we have the differential identities
\beq \label{2.7} 
\frac{d}{dt} \big( \frac{1}{2}| \D u \big(\ga(t,x,\xi) \big)|^{2} \big) = \Big(  \frac{|\D u|^{2}}{| \xi ^ {\top}\D u|^{2}} \xi ^ {\top}  \D u  \ot  \D u :  \D^2 u \Big) \big(\ga(t,x,\xi) \big),
\eeq
\beq \label{2.8} 
\frac{d}{dt} \big(\xi ^ {\top}\D u \big(\ga(t,x,\xi) \big) \big) =| \D u \big(\ga(t,x,\xi) \big)|^{2},
\eeq
which imply $ \D u  \ot  \D u :  \D^2 u = 0$ on $\Om$ if and only if $ |\D u \big(\ga(t,x,\xi) \big)|$ is constant along trajectories $\ga$ and $t \longmapsto \xi ^ {\top}  u \big(\ga(t,x,\xi) \big) $ is affine.
\el

\ms

\ms

\section{Proof of the main result} \label{section3}
In this section we present the proof of the main result of this paper, theorem \ref{theorem1.2}.

\BPT \ref{theorem1.2}. 
  
Let $u \in C^{2}\!\left( \Om,\R^N\right)$ be a solution to the $\infty$-\!Laplace system \eqref{1.1}. Note that the PDE system can be decoupled to the following systems

\beq \label{3.1}
\D u \, \D \Big( \frac{1}{2}|\D u|^{2} \Big)\,=\,0,
\eeq
\beq \label{3.2}
|\D u|^{2} [\![\D u]\!]^\bot \Delta u \,=\,0.
\eeq

Set $\Om_{1} :=\text{ int} \{\rk(\D u)\,\leq\,1\}$,  $\Om_{r} :=\text{ int} \{\rk(\D u)\,\equiv\,r\}$ and $\Om_{n} :=\{\rk(\D u)\,\equiv\,n\}$. Then:
 
On $\Om_{n}$ we have $\rk(\D u)\,=\,\dim(\Om_{n} \sub \R^n)=n$. Since $N\geq n$ and  hence the map $ u\!:\! \Om_{n} \longrightarrow \R^{N}$ is an immersion (because its derivative has constant rank equal to the dimension of the domain, the arguments in the case of $\rk(\D u)\,\equiv\,n$ follows the same lines as in \cite[theorem 1.1]{15} but we provide them for the sake of completeness).  This means that $\D u$ is injective. Thus, $\D u(x)$ possesses a left inverse $(\D u(x))^{-1}$ for all $x \in \Om_{n}$. Therefore, the system \eqref{3.1} implies
\beq \label{3.3}
(\D u)^{-1}\D u \, \D \Big( \frac{1}{2}|\D u|^{2} \Big)\,=\,0,
\eeq
and hence $\D \Big( \frac{1}{2}|\D u|^{2} \Big)\,=\,0$ on $\Om_{n}$, or equivalently
\beq \label{3.4}
|\D u|^{2}\,=\,C^{2},
\eeq
on each connected component of $\Om_{n}$.
Moreover, \eqref{3.4} holds on the common boundary of $\Om_{n}$ with any other component of the partition.

 On $\Om_{r}$  we have $\rk(\D u)\,\equiv\,r$,  where $r$ is integer in $\{2,3,4,...,(n-1)\}$. Consider the gradient flow \eqref{2.6}. Giving that \eqref{3.1} holds, then by the proposition of Gradient flows for tangentially $\infty$ \!\!-\!\! Harmonic maps and its  improved modification lemma which we recalled in the preliminaries, we must have that  $ |\D u \big(\ga(t,x,\xi) \big)|$ is constant along trajectories $\ga$ and $t \longmapsto \xi ^ {\top}  u \big(\ga(t,x,\xi) \big) $ is affine. Moreover, if there exist some $r$ and $q$ such that $2 \leq r < q \leq n-1$, and  $\partial \Om_r \cap \partial \Om_q \, \neq\, \phi$. Then a similar thing happen on  $\partial \Om_r \cap \partial \Om_q \sub S$ (when  both $r\D$ and $q\D$ phases coexist), because in this case we also have that $\rk(\D u)\,\equiv\,r$ and we have same result as above.

\ms
\ms

The proof of the remaining claims of the theorem is very similar to \cite[theorem 1.1]{15}, which we give below for the sake of completeness:

 On $\Om_{1} :=\text{ int} \{\rk(\D u)\,\leq\,1\}$ we have $\rk(\D u)\,\leq\,1$. Hence there exist vector fields $\xi\!:\!  \R^n \supseteq \Om_{1} \longrightarrow \R^{N}$ and  $w\!:\!  \R^n \supseteq \Om_{1} \longrightarrow \R^{n}$ such that $\D u \,=\, \xi \ot w$. Suppose first that $\Om_1$ is contractible. Then, by the Rigidity Theorem  \ref{theorem2.1}, there exist a function $f \in  C^{2}\!\left( \Om_1,\R\right)$, a partition of $\Om_1$ to Borel sets $\big \{ B_i \big \}_{ i\in \N}$ and Lipschitz curves  ${\big\{\mathcal{V}^i \big\} }_{i \in \N} \subseteq W^{1,\infty}_{\text{loc}}(\Om)^N$
 with  $|\dot {\mathcal{V}}^i| \, \equiv \, 1$ on $f(B_i)$, $|\dot {\mathcal{V}}^i| \, \equiv \, 0$ on $\R\!\setminus  \!f(B_i)$ twice differentiable on $ f(B_i)$, such that $u\,=\, \mathcal{V}^i  \circ f $ on each $ B_i  \subseteq \Om$ and hence  $\D u \, = \, ( \mathcal{V}^i  \circ f )   \ot \D f $. By  \eqref{3.1} , we obtain 
\beq \label{3.5} 
\begin{split}
 \Big( \big( \dot { \mathcal{V}^ i} \circ f \big) \ot \D f \Big) \ot & \Big( \big( \dot { \mathcal{V}^ i} \circ f \big) \ot \D f \Big) : \\
& : \Big[  \big( \ddot { \mathcal{V}^ i} \circ f \big) \ot \D f \ot \D f +   \big( \dot { \mathcal{V}^ i} \circ f \big) \ot \D ^2 f \Big]\,=\,0,
\end{split}
\eeq
on $B_i \sub \Om_1$. Since $| \dot { \mathcal{V}^ i}| \equiv 1$ on $f(B_i)$, we have that $ \ddot { \mathcal{V}^ i}$ is normal to $ \dot { \mathcal{V}^ i}$ and hence 
\beq \label{3.6} 
\Big( \big( \dot { \mathcal{V}^ i} \circ f \big) \ot \D f \Big) \ot \Big( \big( \dot { \mathcal{V}^ i} \circ f \big) \ot \D f \Big) : \Big(  \big( \dot { \mathcal{V}^ i} \circ f \big) \ot \D ^2 f \Big)\,=\,0,
\eeq
on $B_i \sub \Om_1$. Hence, by using again that $| \dot { \mathcal{V}^ i}|^2 \equiv 1$ on $f(B_i)$ we get 
\beq \label{3.7} 
\Big(  \D f  \ot \D f  : \D ^2 f \Big) \big( \dot { \mathcal{V}^ i} \circ f \big)\,=\,0,
\eeq
on $B_i \sub \Om_1$. Thus, $ \Delta_\infty f \,=\,0$ on $B_i$. By \eqref{3.2} and again since $| \dot { \mathcal{V}^ i}|^2 \equiv 1$ on $f(B_i)$, we have $ [\![\D u]\!]^\bot \,=\,  [\![ \dot { \mathcal{V}^ i} \circ f]\!]^\bot$ and hence 
\beq \label{3.8} 
|\D f|^2 [\![ \dot { \mathcal{V}^ i} \circ f]\!]^\bot \Div \Big( \big( \dot { \mathcal{V}^ i} \circ f \big) \ot \D f \Big)\,=\,0,  
\eeq
on $B_i \sub \Om_1$. Hence,
\beq \label{3.9} 
|\D f|^2 [\![ \dot { \mathcal{V}^ i} \circ f]\!]^\bot \Big( \big( \ddot { \mathcal{V}^ i} \circ f \big) |\D f|^{2} +  \big( \dot { \mathcal{V}^ i} \circ f \big) \Delta f  \Big)\,=\,0,  
\eeq
on $B_i$, which by using once again  $| \dot { \mathcal{V}^ i}|^2 \equiv 1$  gives
\beq \label{3.10} 
|\D f|^4  \big( \ddot { \mathcal{V}^ i} \circ f \big)\,=\,0,
\eeq
on $B_i$. Since $\Delta_\infty f\,=\,0$ on  $B_i$ and $\Om_1 \,=\, \cup_{1}^{\infty}B_i$, $f$ is $\infty$-Harmonic on $\Om_1$. Thus, by Aronsson's theorem in , either $|\D f| \textgreater 0$ or $|\D f|  \equiv 0$ on $\Om_1$. 

If the first alternative holds, then by \eqref{3.10} we have $\ddot  { \mathcal{V}^ i} \equiv 0$ on $f(B_i)$ for all $i$ and hence, ${ \mathcal{V}^ i}$ is affine on $f(B_i)$, that is  ${ \mathcal{V}^ i}\,=\, t\xi ^i +a^i$ for some $|\xi^i|\,=\,1,  a^i \in \R ^N$. Thus, since $u\,=\, { \mathcal{V}^ i} \circ f $ and $u \in  C^{2}\!\left( \Om_1,\R^N\right) $, all $\xi ^i$ and all $a^i$ coincide and consequently  $u\,=\, \xi f + a$ for $\xi \in \mS ^{N-1}, a \in \R ^N$ and $f \in C^{2}\!\left( \Om_1,\R\right)$.

If the second alternative holds, then $f$ is constant on $\Om_1$ and hence, by the representation $u\,=\,{ \mathcal{V}^ i} \circ f$, $u$ is piecewise constant on each $B_i$. Since $u \in  C^{2}\!\left( \Om_1,\R^N\right)$ and $\Om_1 \,=\, \cup_{i}^{\infty} B_i$, necessarily $u$ is constant on $\Om_1$. But then $|\D u|_{\Om_2}|\,=\,|\D f|_{\mathscr{S}}|\,=\,0$ and necessarily $\Om_2 \,=\, \phi$. Hence, $|\D u|  \equiv 0$ on $\Om$, that is  $u$ is affine on each of the connected components of $\Om$.

If $\Om_1$is not contractible, cover it with balls ${\{\mB _m\}}_{m \in \N}$ and apply the previous argument. Hence, on each $\mB _m$, we have $u\,=\, \xi^m f^m+a^m , \xi ^m \in \mS ^{N-1}, a^m \in \R ^N$ and $f^m \in  C^{2}\!\left( \mB_m,\R\right)$ with $\Delta_\infty f^m\,=\,0$ on  $\mB_m$ and hence either $|\D f^m| > 0$ or $|\D f^m| \equiv 0$. Since $C^{2}\!\left( \Om_1,\R^N\right)$, on the other overlaps of the balls the different expressions of $u$ must coincide and hence, we obtain $u\,=\, \xi f + a$ for $\xi \in \mS ^{N-1}, a \in \R ^N$ and $f \in C^{2}\!\left( \Om_1,\R\right)$ where $\xi$ and $a$ may vary on different connected components of $\Om_1$. The theorem follows.


\qed

\ms
\ms

\subsection*{Acknowledgement.} The author is indebted to N. Katzourakis and would like to thank him for his encouragement and support towards him and for giving the permission to use his own figures, as well as for all the comments and suggestions which improved the content of this paper.

\ms
\ms

\bibliographystyle{abbrv}

\end{document}